\documentclass[a4pper,11pt]{article}
\usepackage{amssymb}
\usepackage{graphicx}
\usepackage{epstopdf}
\usepackage{float}
\usepackage{abstract}
\usepackage{latexsym,bm}
\usepackage{amsfonts}
\usepackage{latexsym}
\usepackage{mathrsfs}
\usepackage{amsmath}
\usepackage{flafter}
\usepackage{color}
\usepackage{subfigure}
\usepackage{geometry}
\usepackage{titlesec}
\usepackage{xcolor}
\usepackage{colortbl,booktabs}  
\usepackage{multicol}  
\usepackage{multirow} 
\usepackage{appendix}
\usepackage{indentfirst}
\usepackage[numbers,sort,compress]{natbib}
\usepackage[utf8]{inputenc}

\numberwithin{equation}{section}

\newcommand{\bea}{\begin{eqnarray*}}
\newcommand{\eea}{\end{eqnarray*}}
\newcommand{\be}{\begin{equation}}
\newcommand{\ee}{\end{equation}}
\newcommand{\ben}{\begin{eqnarray*}}
\newcommand{\een}{\end{eqnarray*}}

\begin{document}
\title { Optimal error estimate of a conservative Fourier pseudo-spectral method for the  space fractional nonlinear Schr\"{o}dinger equation}
\author{ Zhuangzhi Xu$^1$, Wenjun Cai$^1$, Chaolong Jiang$^2$, Yushun Wang$^1$\footnote{Correspondence author. Email:wangyushun@njnu.edu.cn.}
\\{\small$^1$ Jiangsu Key Laboratory for NSLSCS,}
\\{\small School of Mathematical Sciences,  Nanjing Normal University,}
\\{\small  Nanjing 210023, China}\\
{\small $^2$ School of Statistics and Mathematics, }\\
{\small Yunnan University of Finance and Economics, Kunming 650221, P.R. China}\\}
\date{}
\maketitle

\begin{abstract}

 In this paper, we consider the error analysis of a conservative Fourier pseudo-spectral method that conserves mass and energy for the space fractional nonlinear Schr\"{o}dinger equation. We give a new fractional Sobolev norm that can construct the discrete fractional Sobolev space, and we also can prove some important lemmas for the new fractional Sobolev norm. Based on these lemmas and energy method,
  a priori error estimate for the method can be established.
 Then, we are able to prove that 
 the Fourier pseudo-spectral method is unconditionally convergent with order $O(\tau^{2}+N^{\alpha/2-r})$ in the  discrete $L^{\infty}$ norm,
where $\tau$ is the time step and $N$ is the number of collocation points used in the spectral method. Numerical examples are presented to verify the theoretical analysis.
\end{abstract}
\textbf{AMS subject classification:} 35R11, 65M70\\[2ex]
\textbf{Keywords:} Fractional nonlinear Schr\"{o}dinger equation; Consevation laws;
 Fourier pseudo-spectral method; Priori error estimate

\section{Introduction}
% Fractional Schr\"{o}dinger equation is the fractional generalization of Schr\"{o}dinger equation. It was first introduced by Laskin~\cite{laskin2000fractionalL,laskin2000fractionalR} via defining fractional path integral approach over the path of the L$\acute{\rm e}$vy flights instead of Feynman path in the development of fractional quantum and statistical mechanics. 

% There are some interesting applicaitions of fractional  Schr\"{o}dinger equation in physical systems, including nonlinear optics~\cite{longhi2015fractional}, quantum physics~\cite{guo2006some}, Propagation dynamics~\cite{zhang2016propagation} and water wave dynamics~\cite{ionescu2014nonlinear}. Also, from mathematical viewpoints, most of the studies focus on the well-posedness of fractional Schrodinger equation~\cite{guo2013well,cheng2012bound,guo2008existence,hu2011global}.  
The nonlinear fractional Schr\"{o}dinger equation is a generalization of the classical Schr\"{o}dinger equation. It has found several applicaitions in physics, such as nonlinear optics~\cite{longhi2015fractional}, propagation dynamics~\cite{zhang2016propagation} and water wave dynamics~\cite{ionescu2014nonlinear}.
In this paper,  we consider the following space fractional nonlinear Schr\"{o}dinger~(FNLS) equation 
   \begin{align}\label{aa1}
    \text{i} u_{t}-(-\Delta)^{\frac{\alpha}{2}} u+\beta|u|^{2}u=0,~~x\in \Omega,~~0<t\leq T,
   \end{align}
with the periodic boundary condition
\begin{equation}
   \begin{aligned}\label{aa2}
   u(x,t)=u(x+L,t),\quad x\in \Omega,~0<t\leq T,
  \end{aligned}
\end{equation}
and the initial condition
   \begin{align}\label{aa3}
   u(x,0)=\varphi(x),~~ x\in \Omega,
   \end{align}
where $\text{i}=\sqrt{-1},1<\alpha\leq 2$,  $\Omega=\left[a,b\right]$ and $L=b-a.$
 $u(x,t)$ is a complex-valued wave function, parameter $\beta$ is a real constant,
 and $\varphi(x)$ is a complex-value initial data. 
The fractional Laplacian $(-\triangle)^{\frac{\alpha}{2}}$ acting on periodic function defined by~\cite{pozrikidis2016fractional}
\begin{align}\label{aa4}
(-\triangle)^{\frac{\alpha}{2}}u=\sum_{k\in \mathbb{Z}}|\mu k|^{\alpha}\widehat{u}_{k}e^{\text{i}\mu kx},~~ \mu=\frac{2\pi}{L},
\end{align}
where 
\begin{align}\label{aa5}
  u=\sum_{k\in \mathbb{Z}}\widehat{u}_{k}e^{\text{i}\mu kx},~~ \widehat{u}_{k}=\frac{1}{L}\int_{\Omega}u(x)e^{-\text{i}\mu kx}dx.
\end{align}

 When $\alpha=2$, the equation~(\ref{aa1}) reduces to the classical nonlinear Schr\"{o}dinger~(NLS) equation. Due to self-adjoint property of the fractional Laplacian, the solution of $(\ref{aa1})\verb|-|(\ref{aa3})$ satisfies the following mass and energy conservation laws:
\begin{align}\label{aa6}
{\rm Mass}:~~ M(t)=\int_{\Omega}|u(x,t)|^{2}dx=M(0),~~~~~~~~~~\quad\quad\quad\quad\quad\\
{\rm Energy}:~~ E(t)=\int_{\Omega}|(-\Delta)^{\frac{\alpha}{4}}u(x,t)|^{2}-\frac{\beta}{2}|u(x,t)|^{4}dx=E(0).\label{aa7}
\end{align}

  	Various numerical methods have been developed in the literatures for the space FNLS equation, including finite difference methods~\cite{wang2014linearly,wang2013crank,zhao2014fourth,wang2016point,wang2015maximum}, finite element methods~\cite{li2018fast,li2017unconditionally}, spectral methods~\cite{wang2019split,bhrawy2015fully}.
  In the past few decades, structure-preserving methods which can inherit the intrinsic geometric properties of the given dynamical system have attracted a lot of interest due to the superior properties in long time numerical simulation over traditional methods. For more details, readers can refer  to~\cite{feng2010symplectic,hairer2006geometric,leimkuhler2004simulating}.
  Recently, structure-preserving numerical methods have been extended to solve the space FNLS equation. For example, in \cite{wang2013crank},
   Wang et al. first constructed a mass conservative Crank-Nicolson difference scheme, and 
   they further proposed a linearly implicit difference scheme that conserves mass and energy in~\cite{wang2014linearly}. In~\cite{wang2015maximum,wang2016point},
   the modified mass and energy conservative Crank-Nicolson difference schemes were presented.  Other works related to the conservative method can be founded in~\cite{ran2016conservative,li2018fast}.

Spectral and pseudo-spectral methods have been proved to be an efficient and high order numerical method in solving smooth problems~\cite{tang2006spectral}. Over the past few years,
though structure-preserving spectral methods have been widely used to solve Hamiltonian PDEs~\cite{Bridges2001Multi,chen2001multi,kong2010semi,gong2014multi}.
Only in very recently years, structure-preserving Fourier pseudo-spectral methods are extended to solve the space FNLS equation. For instance, in \cite{wang2018structure},  Wang and Huang proposed
 the symplectic and multi-symplectic Fourier pseudo-spectral methods for the space FNLS equation.
 In~\cite{miyatake2019structure}, a mass and energy conservative Fourier pseudo-spectral method was constructed. The numerical results show that the conservative Fourier pseudo-spectral method is efficient and stable for long-term numerical simulation. However the unconditionally convergent results on the conservative Fourier pseudo-spectral method for the space fractional PDEs have not been obtained.  Actually, with the help of the defined fractional Sobolev norm and the discrete uniform Gagliardo-Nirenberg inequality in~\cite{kirkpatrick2013continuum}, we can easily prove that the conservative Fourier pseudo-spectral method for the space FNLS equation is unconditionally convergent in the discrete $L^{2}$ norm, but the challange problem is error estimate in  $L^{\infty}$ norm. For the classical NLS equations,
 in~\cite{gong2017conservative}, Gong first established the semi-norm equivalence between the finite difference method and the Fourier pseudo-spectral method and thus obtained the unconditionally convergent results on the Fourier pseudo-spectral in the discrete  $L^{2}$ norm. Then based on this equivalence, the error estimates of the Fourier pseudo-spectral method in the discrete $L^{\infty}$ norm 
were obtained in~\cite{jiang2018optimal}.
However, this error analysis technique for establishing semi-norm equivalence
can not extend to the FNLS equations.
By reading the finite difference methods for the fractional PDEs with fractional Laplacian \cite{zhao2014fourth,wang2016implicit,hou2017numerical}, we know that under the homogeneous Dirichlet boundary condition, the
 Riesz derivative is discreted instead of fractional Laplacian due to the equivalence between Riesz derivative and fractional Laplacian.
But this equivalence does not hold in the case of periodic boundary condition. That's why even now these is no corresponding finite difference method has been used to solve the space fractional PDEs under the periodic boundary condition. Therefore, for the FNLS equation, it is impossible to establish semi-norm equivalence between the finite difference method and the Fourier pseudo-spectral method in the error analysis.

To obtain the  $L^{\infty}$ norm error estimates of the Fourier pseudo-spectral method for the space FNLS equation, in this paper, we introduce
the discrete fractional Sobolev space $H^{\alpha/2}_{h}$ with a new  discrete fractional Sobolev norm. We establish several lemmas for the new discrete fractional norm, 
 based on these important lemmas and the energy method, a prior estimate for the method is estimated. Then we can prove that the conservative Fourier pseudo-spectral method is unconditionally convergent with order of $O(\tau^{2}+N^{\alpha/2-r})$ in the  disctete $L^{\infty}$ norm.

The rest of the paper is organized as following. In section 2, we construct the discrete fractional Sobolev space by introducing a new fractional Sobolev norm and we also prove some important lemmas for the new fractional Sobolev norm.
In section 3,  a conservative Fourier pseudo-spectral scheme for the FNLS equation is given, we show the numerical scheme satisfies discrete conservation laws and obtain a priori estimate. In section 4, the convergence property of the scheme is analyzed. Subsequently in section 5, we carry out some numerical experiments to confirm our theoretical results and show the efficiency of the scheme. Finally, we give a conclusion  in section 6.
\section{Fourier pseudo-spectral method}
\newtheorem{theorem}{Theorem}[section]
\newtheorem{remark}{Remark}[section]
Let $N$ be an even integer, we define  step size in space: $h=L/N$. Then, the spatial grid points are defined as follows:$~\Omega_{h}=\{x_{j}=a+jh,j=0,1,\ldots, N\}$. For any positive integer $N_{t}$, we define the time-step: $\tau=T/N_{t}$. Then grid points in space and time are given by $\Omega_{h\tau}=\Omega_{h}\times\Omega_{\tau}$, where $\Omega_{\tau}=\left\{t_{n}=n\tau,n=0,1,\ldots,N_{t}\right\}$.
For a grid function $u=\{u^{n}_{j}|(x_{j},t_{n})\in\Omega_{h\tau}\}$, we introduce the following notations:
 $$\delta^{+}_{x}u^{n}_{j}=\frac{u^{n}_{j+1}-u^{n}_{j}}{h},\quad u^{n+\frac{1}{2}}_{j}=\frac{u^{n+1}_{j}+u^{n}_{j}}{2},\quad \delta^{+}_{t}u^{n}_{j}=\frac{u^{n+1}_{j}-u^{n}_{j}}{\tau}.$$
 
 Let $\mathcal{V}_{h}=\{u|u=(u_{j}),x_{j}\in \Omega_{h}\}$ be the space of grid functions defined on $\Omega_{h}$. For any grid function $u,v\in \mathcal{V}_{h}$, we define the discrete inner product and associated $L^{2}$ norm
 \begin{align}\label{eq1}
 (u,v)_{h}=\dfrac{1}{N}\sum\limits_{j=0}^{N-1}u_{j}\overline{v}_{j},~~\|u\|^{2}_{h}=(u,u)_{h}.
 \end{align}
 We also define the discrete $L^{p}$ norm as
 \begin{align}
 \|u\|^{p}_{l^{p}_{h}}=\dfrac{1}{N}\sum\limits_{j=0}^{N-1}|u_{j}|^{p},~~~1\leq p<+\infty,
 \end{align}
 and the discrete $L^{\infty}$ norm as
 \begin{align}
 \|u\|_{l^{\infty}_{h}}=\max\limits_{0\leq j \leq N-1}|u_{j}|.
 \end{align}
 % so that $I_{N}u(x_{j})=u_{j},~$ where $u_{j}=u(x_{j}),~~g_{j}(x_{l})=\delta^{j}_{l}.$
 
 %For the sake of consistency, we consider the Fourier interpolation on $N$ collocation points $\{x_{j}\}^{N-1}_{j=0}$, denote the interpolation operator by $I_{N}$. That is,
 %\begin{align}
 % (I_{N}u)(x)=\sum_{k=-N/2}^{N/2-1}\widehat{u}_{k}e^{ik\mu x},\quad \mu=\frac{2\pi}{L},
 %\end{align}
 \subsection{Disctete Fractional Sobolev norm} 
 We define a function space $S_{N}$ by
 $$S_{N}=\text{span}\{g_{j}(x),~~j=0,1,\cdots,N-1\},$$
 where $g_{j}(x)$ is a trigonometric polynomial defined by
\begin{align}
 g_{j}(x)=\dfrac{1}{N}\sum_{k=-N/2}^{N/2}\frac{1}{c_{k}}e^{\text{i}k\mu(x-x_{j})},
 \end{align}
 where
$$
 c_{k}=\left\{\begin{aligned}
 &1,~|k|<N/2,\\
 &2,~|k|=N/2,
 \end{aligned}~~~~\mu=\frac{2\pi}{L}.
 \right.
$$
 Then, we define the interpolation operator $I_{N}:L^{2}(\Omega)\rightarrow S_{N}$ by
 \begin{equation}
 \begin{aligned}\label{a1}
 I_{N}u(x)=\sum_{j=0}^{N-1}u_{j}g_{j}(x)=\sum_{k=-N/2}^{N/2}\widehat{u}_{k}e^{\text{i}k\mu x},
 \end{aligned}
 \end{equation}
where
  \begin{align}\label{b1}
    \widehat{u}_{k}=\dfrac{1}{Nc_{k}}\sum\limits_{j=0}^{N-1}u_{j}e^{-\text{i}k\mu x_{j}},~~-N/2 \leq k \leq N/2,
    \end{align}
and $\widehat{u}_{\frac{N}{2}}=\widehat{u}_{-\frac{N}{2}}$ for $k=\frac{N}{2}$. Therefore we have the inverse transformation
\begin{align}\label{b2}
u_{j}=(I_{N}u)(x_{j})=\sum_{k=-N/2}^{N/2-1}\widehat{u}_{k}e^{\text{i}k\mu x_{j}}.
\end{align}
%According to $(\ref{eq1})$ and $(\ref{b1})-(\ref{b2})$, 
For any $u\in l^{2}_{h}:=\{u|u\in\mathcal{V}_{h},\|u\|^{2}_{h}<\infty\}$, we have $\widehat{u}\in l^{2}:=\{x=\{x_{k}\}|\sum\limits_{k=-\infty}^{\infty}x^{2}_{k}<\infty\},$ and the Parseval's theorem gives
\begin{align}\label{p}
  (u,v)_{h}=\sum\limits_{k=-N/2}^{N/2-1}\widehat{u}_{k}\overline{\widehat{v}}_{k}.
\end{align}

Given a constant $\sigma\in[0,1]$, we define the discrete fractional Sobolev norm $\|\cdot\|_{H^{\sigma}_{h}}$ and semi-norm $|\cdot|_{H^{\sigma}_{h}} $ as
\begin{align}\label{2}
  |u|^{2}_{H^{\sigma}_{h}}=\sum\limits_{k=-N/2}^{N/2-1}|\mu k|^{2\sigma}|\widehat{u}_{k}|^{2},~~~~~~~~~~~~~~ \|u\|^{2}_{H^{\sigma}_{h}}=\sum\limits_{k=-N/2}^{N/2-1}(1+|\mu k|^{2\sigma})|\widehat{u}_{k}|^{2}.
\end{align}
 Clearly, $\|u\|^{2}_{H^{\sigma}_{h}}=\|u\|^{2}_{h}+|u|^{2}_{H^{\sigma}_{h}},\|u\|^{2}_{H^{0}_{h}}=\|u\|^{2}_{h}.$
 We can easily prove that the discrete Sobolev spaces is the normed linear spaces according to the norm $\|u\|_{H^{\sigma}_{h}}$ defined in $(\ref{2})$.
Next, we introduce the following lemmas, which are important for unconditional convergence analysis of the conservative Fourier pseudo-spectral method.
\newtheorem{lemma}{Lemma}[section]
\begin{lemma} \label{lemma1}
{\rm(Discrete uniform Sobolev inequality).} For any $\frac{1}{2}<\sigma\leq 1,$ there exists a constant $C=C(\sigma)>0$ independent of $h>0$ such that
\begin{equation}\label{lem}
  \|u\|_{l^{\infty}_{h}}\leq C\|u\|_{H^{\sigma}_{h}}.
\end{equation}
\end{lemma}
{\bf Proof.}\quad From the inverse transformation $(\ref{b2})$ and the Cauchy-Schwarz inequality, we obtain
\begin{equation}
  \begin{aligned}\label{1}
     \|u\|_{l^{\infty}_{h}}&\leq \sum_{k=-N/2}^{N/2-1}|\widehat{u}_{k}|\\
     &=\sum_{k=-N/2}^{N/2-1}\frac{1}{(1+|\mu k|^{2\sigma})^{\frac{1}{2}}}(1+|\mu k|^{2\sigma})^{\frac{1}{2}}|\widehat{u}_{k}|\\
     &\leq \bigg(\sum_{k=-N/2}^{N/2-1}\frac{1}{1+|\mu k|^{2\sigma}}\bigg)^{\frac{1}{2}}\bigg(\sum_{k=- N/2}^{N/2-1}(1+|\mu k|^{2\sigma})|\widehat{u}_{k}|^{2}\bigg)^{\frac{1}{2}}\\
     &\leq \bigg(\sum_{k=-N/2}^{N/2-1}\frac{1}{1+|\mu k|^{2\sigma}}\bigg)^{\frac{1}{2}}\|u\|_{H^{\sigma}_{h}}.
   \end{aligned}
\end{equation}
For $\frac{1}{2}<\sigma\leq 1,$ this implies $(\ref{lem})$ and thus the proof is completed.
\begin{lemma}\label{lemmaa2}
 For $0\leq \sigma_{0}\leq \sigma \leq 1,$ there exist a constant $C\in[1,2]$ such that
\begin{align}\label{lem1}
  \|u\|_{H^{\sigma_{0}}_{h}}\leq C\|u\|^{\frac{\sigma_{0}}{\sigma}}_{H^{\sigma}_{h}}\|u\|^{1-\frac{\sigma_{0}}{\sigma}}_{h}.
\end{align}
\end{lemma}
{\bf Proof.}\quad From the definition of $\|u\|_{H^{\sigma_{0}}_{h}}$ and the H\"{o}lder's inequality, we have
\begin{equation}
  \begin{aligned}\label{3}
  \|u\|^{2}_{H^{\sigma_{0}}_{h}}&=\sum_{k=-N/2}^{N/2-1}(1+|\mu k|^{2\sigma_{0}})|\widehat{u}_{k}|^{2}\\
  &=\sum_{k=-N/2}^{N/2-1}\bigg((1+|\mu k|^{2\sigma})|\widehat{u}_{k}|^{2}\bigg)^{\frac{\sigma_{0}}{\sigma}}(|\widehat{u}_{k}|^{2})^{1-\frac{\sigma_{0}}{\sigma}}\bigg(\frac{1+|\mu k|^{2\sigma_{0}}}{(1+|\mu k|^{2\sigma})^{\frac{\sigma_{0}}{\sigma}}}\bigg)\\
  &\leq C \bigg(\sum_{k=-N/2}^{N/2-1}(1+|\mu k|^{2\sigma})|\widehat{u}_{k}|^{2}\bigg)^{\frac{\sigma_{0}}{\sigma}}\bigg(\sum_{k=-N/2}^{N/2-1}|\widehat{u}_{k}|^{2}\bigg)^{1-\frac{\sigma_{0}}{\sigma}}\\
  &=C\bigg(\sum_{k=-N/2}^{N/2-1}(1+|\mu k|^{2\sigma})|\widehat{u}_{k}|^{2}\bigg)^{\frac{\sigma_{0}}{\sigma}}\bigg(\sum_{k=-N/2}^{N/2-1}|\widehat{u}_{k}|^{2}\bigg)^{1-\frac{\sigma_{0}}{\sigma}}\\
   &=C\bigg(\|u\|^{\frac{\sigma_{0}}{\sigma}}_{H^{\sigma}_{h}}\|u\|^{1-\frac{\sigma_{0}}{\sigma}}_{h}\bigg)^{2},
  \end{aligned}
\end{equation}
where  the inequality holds due to the fact $\frac{1}{2}(1+a^{\mu})\leq (1+a)^{\mu}\leq (1+a^{\mu})$ for $a>0,~0\leq \mu \leq 1$. Thus the proof is completed.
\begin{lemma}{\rm( Hausdorff-Young inequality)}.\label{lemma6}
	If $1\leq q\leq 2,~\frac{1}{q}+\frac{1}{p}=1,$ then
	\begin{align}
	\bigg(h\sum_{j=0}^{N-1}|u_{j}|^{p}\bigg)^{\frac{1}{p}}\leq \bigg(\sum_{k=-N/2}^{N/2-1}|\widehat{u}_{k}|^{q}\bigg)^{\frac{1}{q}}.
	\end{align}
\end{lemma}
{\bf Proof.}\quad From the inverse transformation $(\ref{b2})$, we have
\begin{equation}
\sup\limits_{0\leq j\leq N-1}|u_{j}|\leq \sum\limits_{k=-N/2}^{N/2-1}|\widehat{u}_{k}|,
\end{equation}
the Parseval's identity gives
\begin{equation}
h\sum_{j=0}^{N-1}|u_{j}|^{2}=\sum_{k=-N/2}^{N/2-1}|\widehat{u}_{k}|^{2}.
\end{equation}
Then using the Riesz-Thorin Interpolation theorem~(see Theorem 8.6 in [\citealp[page 316]{castillo2016introductory}]), we can obtain the conclusion.
\begin{lemma} \label{lemma3}
 For any $\frac{p-2}{2p}<\sigma_{0}\leq 1,$ there exists a constant $C_{\sigma_{0}}=C(\sigma_{0})>0$ independent of $h>0,$ such that
\begin{equation}\label{lem2}
  \|u\|_{l^{p}_{h}}\leq C_{\sigma_{0}}\|u\|^{\frac{\sigma_{0}}{\sigma}}_{H^{\sigma}_{h}}\|u\|^{1-\frac{\sigma_{0}}{\sigma}}_{h},~~~~2\leq p\leq +\infty,~~\sigma_{0}\leq \sigma \leq 1.
\end{equation}
\end{lemma}
{\bf Proof.}\quad By Lemma $\rm\ref{lemma6}$ and H\"{o}lder's inequality, for $1\leq q \leq 2$ such that $\frac{1}{p}+\frac{1}{q}=1,$ we have
\begin{equation}
  \begin{aligned}\label{4}
 \bigg(h\sum_{j=0}^{N-1}|u_{j}|^{p}\bigg)^{\frac{1}{p}}&\leq \bigg(\sum_{k=-N/2}^{N/2-1}|\widehat{u}_{k}|^{q}\bigg)^{\frac{1}{q}}\\
 &=\bigg(\sum_{k=-N/2}^{N/2-1}\frac{1}{(1+|\mu k|^{2\sigma_{0}})^{\frac{q}{2}}}(1+|\mu k|^{2\sigma_{0}})^{\frac{q}{2}}|\widehat{u}_{k}|^{q}\bigg)^{\frac{1}{q}}\\
 &\leq \bigg(\sum_{k=-N/2}^{N/2-1}(1+|\mu k|^{2\sigma_{0}})|\widehat{u}_{k}|^{2}\bigg)^{\frac{1}{2}}\bigg(\sum_{k=-N/2}^{N/2-1}\frac{1}{(1+|\mu k|^{2\sigma_{0}})^{\frac{q}{2-q}}}\bigg)^{\frac{2-q}{2q}}\\
 & \leq \|u\|_{H^{\sigma_{0}}_{h}}\bigg(\sum_{k=-N/2}^{N/2-1}\frac{1}{(1+|\mu k|^{2\sigma_{0}})^{\frac{q}{2-q}}}\bigg)^{\frac{2-q}{2q}}.
  \end{aligned}
\end{equation}
Then for $\frac{p-2}{2p}<\sigma_{0}\leq 1,$ we have
\begin{align}
\|u\|_{l^{p}_{h}}\leq \tilde{C}_{\sigma_{0}}\|u\|_{H^{\sigma_{0}}_{h}},
\end{align}
where $\tilde{C}_{\sigma_{0}}=\tilde{C}(\sigma_{0})>0$ is independent of $h$. Combining the above inequality with $(\ref{lem1})$ gives $(\ref{lem2})$ and thus completes the proof.

%\begin{remark}\rm
%The proof processes of Lemma $\ref{lemma1}$ and Lemma $\ref{lemmaa2}$, Lemma $\ref{lemma3}$ are similar to that in [\citealp[Lemma 3.1]{kirkpatrick2013continuum}] and [\citealp[Lemma 2.2, Lemma 2.3 ]{wang2016implicit}], respectively.
%\end{remark}
\subsection{Discrete fractional Laplacian}
Applying the fractional Laplacian $(-\Delta)^{\frac{\alpha}{2}}$ to the interpolated function $(\ref{a1})$ yields
\begin{equation}\label{a2}
  \begin{aligned}
  (-\Delta)^{\frac{\alpha}{2}}I_{N}u(x)=\frac{1}{N}\sum_{j=0}^{N-1}u_{j}\sum_{p=-N/2}^{N/2}\frac{1}{c_{p}}|\mu p|^{\alpha}e^{\text{i}p\mu(x-x_{j})},
  \end{aligned}
\end{equation}
and thus
\begin{equation}\label{a3}
  \begin{aligned}
  (-\Delta)^{\frac{\alpha}{2}}I_{N}u(x_{k})&=\sum_{p=-N/2}^{N/2-1}d_{p}\bigg(\frac{1}{N}\sum_{j=0}^{N-1}u_{j}e^{-\frac{2\pi \text{i}jp}{N}}\bigg)e^{\frac{2\pi \text{i}pk}{N}},\\
  \end{aligned}
\end{equation}
where
\begin{equation}
  d_{p}=|\mu p|^{\alpha},\quad\quad -N/2\leq p\leq N/2-1.
\end{equation}
%$$d_{p}=\left\{\begin{aligned}
%&|\mu p|^{\alpha},\quad 0\leq p\leq \frac{N}{2},\\
%&|\mu (p-N)|^{\alpha},\quad \frac{N}{2}\leq p\leq N-1.
%\end{aligned}
%\right.$$
For $U\in \mathcal{V}_{h}$,~~we define a discrete fractional Laplacian $(-\Delta)^{\frac{\alpha}{2}}_{d}$ by
\begin{equation}\label{a4}
  \begin{aligned}
  ((-\Delta)^{\frac{\alpha}{2}}_{d}U)_{k}&=\sum_{p=-N/2}^{N/2-1}d_{p}\bigg(\frac{1}{N}\sum_{j=0}^{N-1}U_{j}e^{-\frac{2\pi \text{i}jp}{N}}\bigg)e^{\frac{2\pi \text{i}pk}{N}},\\
  \end{aligned}
\end{equation}
By using the notation of the discrete Fourier transform and its inverse:
\begin{equation}\label{a5}
  \begin{aligned}
  (\mathcal{F}_{d}U)_{k}=\frac{1}{N}\sum_{j=0}^{N-1}U_{j}e^{-\frac{2\pi \text{i}jk}{N}},\quad (\mathcal{F}^{-1}_{d}\widehat{U})_{j}=\sum_{k=-N/2}^{N/2-1}\widehat{U}_{k}e^{\frac{2\pi \text{i}jk}{N}},
  \end{aligned}
\end{equation}
the discrete fractional Laplacian can be expressed as
\begin{equation}\label{a6}
  \begin{aligned}
  (-\Delta)^{\frac{\alpha}{2}}_{d}U=\mathcal{F}^{-1}_{d}{\Lambda_{\alpha}}\mathcal{F}_{d}U,
  \end{aligned}
\end{equation}
where ${\Lambda_{\alpha}}=diag(d_{-\frac{N}{2}},d_{-(\frac{N}{2}-1)},\ldots,0,1,\ldots,d_{\frac{N}{2}-1}).$
Next, we  give several lemmas that show the relationship between discrete fractional Soboolv semi-norm and fractional Laplacian.
%\bm{\wedge_{\alpha}}=diag(|\mu \frac{N}{2}|^{\alpha},|\mu (\frac{N}{2}-1)|^{\alpha},\ldots,0,1,\ldots,|\mu (\frac{N}{2}-1)|^{\alpha}). 
\begin{lemma}\label{lemu}
	For any grid function $u\in \mathcal{V}_{h}$, we have
	\begin{equation}\label{2.03}
	(D_{\alpha}u,u)_{h}=|u|^{2}_{H^{\alpha/2}_{h}}, \quad 1<\alpha\leq 2.
	\end{equation}
\end{lemma}
{\bf Proof.}\quad Using the Parseval's identity $(\ref{p})$, we have
\begin{equation}\label{a8}
\begin{aligned}
(D_{\alpha}u,u)_{h}&=(\mathcal{F}^{-1}_{d}{\Lambda_{\alpha}}\mathcal{F}_{d}u,u)_{h}\\
&=\sum_{k=-N/2}^{N/2-1}({\Lambda_{\alpha}}\mathcal{F}_{d}u)_{k}(\mathcal{F}_{d}\overline{u})_{k}
=\sum_{k=-N/2}^{N/2-1}d_{k}\widehat{u}_{k}\overline{\widehat{u}}_{k}=|u|^{2}_{H^{\alpha/2}_{h}}.
\end{aligned}
\end{equation}
\begin{lemma}\label{2.04}
	For any two grid functions $u,v\in \mathcal{V}_{h}$, we have
	\begin{equation}\label{2.05}
	(D_{\alpha}u,v)_{h}=(D_{\alpha/2}u,D_{\alpha/2}v)_{h},\quad 1<\alpha\leq2.
	\end{equation}
\end{lemma}
{\bf Proof.}\quad Using the Parseval's identity $(\ref{p})$, we have
\begin{equation}\label{a9}
\begin{aligned}
(D_{\alpha}u,v)_{h}&=
(\mathcal{F}^{-1}_{d}{\Lambda_{\alpha}}\mathcal{F}_{d}u,v)_{h}\\
&=\sum_{k=-N/2}^{N/2-1}({\Lambda_{\alpha}}\mathcal{F}_{d}u)_{k}(\mathcal{F}_{d}\overline{v})_{k}\\
&= \sum_{k=-N/2}^{N/2-1}({\Lambda_{\frac{\alpha}{2}}}\mathcal{F}_{d}u)_{k}({\Lambda_{\frac{\alpha}{2}}}\mathcal{F}_{d}\overline{v})_{k}
=(D_{\alpha/2}u,D_{\alpha/2}v)_{h}.
\end{aligned}
\end{equation}
\begin{lemma}\label{2.06}
	For any two grid functions $u,v\in \mathcal{V}_{h}$, we have
	\begin{equation}\label{2.07}
	(D_{\alpha}u,v)_{h}\leq |u|_{H^{\alpha/2}_{h}}|v|_{H^{\alpha/2}_{h}},\quad 1<\alpha\leq2.
	\end{equation}
\end{lemma}
{\bf Proof.}\quad Using the Parseval's identity $(\ref{p})$, we have
\begin{equation}\label{a10}
\begin{aligned}
(D_{\alpha}u,v)_{h}&=(\mathcal{F}^{-1}_{d}{\Lambda_{\alpha}}\mathcal{F}_{d}u,v)_{h}\\
&=\sum_{k=-N/2}^{N/2-1}({\Lambda_{\alpha}}\mathcal{F}_{d}u)_{k}(\mathcal{F}_{d}\overline{v})_{k}
=\sum_{k=-N/2}^{N/2-1}d_{k}\widehat{u}_{k}\overline{\widehat{v}}_{k}.
\end{aligned}
\end{equation}
Therefore
\begin{equation}\label{a11}
\begin{aligned}
(D_{\alpha}u,v)_{h}
% &=\sum_{k=-N/2}^{N/2-1}d_{k}\widehat{u}_{k}\overline{\widehat{v}_{k}}\\
&\leq \bigg(\sum_{k=-N/2}^{N/2-1}d_{k}|\widehat{u}_{k}|^{2}\bigg)^{\frac{1}{2}}\bigg(\sum_{k=-N/2}^{N/2-1}d_{k}|\widehat{v}_{k}|^{2}\bigg)^{\frac{1}{2}}
=|u|_{H^{\alpha/2}_{h}}|v|_{H^{\alpha/2}_{h}}.
\end{aligned}
\end{equation}
 For simplicity, we denote $u^{n}_{j}=u(x_{j},t_{n})$ and $U^{n}_{j}$ as the exact value of $u(x,t)$ and its numerical approximation at $(x_{j},t_{n}),$ respectively.

\section{Solution existence and conservation of the scheme}
We discretize the FNLS equation $(\ref{aa1})\verb|-|(\ref{aa3})$ using the Fourier pseudo-spectral method in space and the Crank-Nicolson method in time to arrive at a fully discrete  system:
\begin{equation}\label{2.01}
  \text{i}\delta^{+}_{t}U^{n}_{j}-(D_{\alpha}U^{n+1/2})_{j}+\frac{\beta}{2}(|U^{n}_{j}|^{2}+|U^{n+1}_{j}|^{2})U^{n+1/2}_{j}=0,\quad U^{n}\in \mathcal{V}_{h},
\end{equation}
where $D_{\alpha}=(-\Delta)^{\frac{\alpha}{2}}_{d}=\mathcal{F}^{-1}_{d}{\Lambda_{\alpha}}\mathcal{F}_{d},~j=0,1,\ldots,N-1.$ For convenience, scheme $(\ref{2.01})$ can be written in an equivalent form
\begin{equation}\label{2.02}
  \text{i}\delta^{+}_{t}U^{n}-D_{\alpha}U^{n+1/2}+F(U^{n},U^{n+1})=0, \quad U^{n}\in\mathcal{V}_{h},
\end{equation}
where $U^{n}=(U^{n}_{j}),~ F(U^{n},U^{n+1})=F(U^{n}_{j},U^{n+1}_{j})=\big(\frac{\beta}{4}(|U^{n}_{j}|^{2}+|U^{n+1}_{j}|^{2})(U^{n}_{j}+U^{n+1}_{j})\big).$
%{\bf Remark.}
%Lemma $\ref{lemu}$ means that we need not establish the equivalence between the difference method and the Fourier pseudo-spectral method as compared with the proof process of error estimate of Fourier pseudo-spectral method for  the classical Schrodinger method(see Lemma 2.4 in \cite{gong2017conservative}).

\begin{lemma}\label{3.01}
 For the approximation $U^{n}\in \mathcal{V}_{h},$ there exist identities:
\begin{gather}
\textrm{\rm{Im}}(D_{\alpha}U^{n+1/2},U^{n+1/2})_{h}=0,~~~~~~~~~~~~~~~~~~~~~~~~~~~~~\\
\textrm{\rm{Re}}(D_{\alpha}U^{n+1/2},\delta^{+}_{t}U^{n})_{h}=\frac{1}{2\tau}(|U^{n+1}|^{2}_{H^{\alpha/2}_{h}}-|U^{n}|^{2}_{H^{\alpha/2}_{h}}),
\end{gather}
\end{lemma}
According to Lemma $\ref{lemu}$ and Lemma $\ref{2.04}$, we can get the results immediately.
Here \textquotedblleft$\text{Im}(s)$\textquotedblright and \textquotedblleft$\rm{ Re}(s)$\textquotedblright  mean taking the imaginary part and real part of a complex number $s$, respectively.
\subsection{Conservation}
\begin{theorem}
The scheme $(\ref{2.02})$ is conservative in the sense that
\begin{gather}\label{3.04}
M^{n}=M^{0},\quad 0\leq n\leq N, \\
E^{n}=E^{0}, \quad 0\leq n\leq N,  \label{3.05}
\end{gather}
\end{theorem}
where
\begin{equation}\label{3.02}
  M^{n}:=\|U^{n}\|^{2}_{h}, \quad\quad E^{n}=|U|^{2}_{H^{\alpha/2}_{h}}-\frac{\beta}{2}\|U^{n}\|^{4}_{l^{4}_{h}}
\end{equation}
{\bf Proof.}\quad Computing the discrete inner product of $(\ref{2.02})$ with $U^{n+1/2}$, then taking the imaginary part, we obtain
\begin{equation}\label{3.03}
  \frac{1}{2\tau}(\|U^{n+1}\|^{2}_{h}-\|U^{n}\|^{2}_{h})=0,\quad  t_{n}\in \Omega_{\tau},
\end{equation}
where Lemma $\ref{3.01}$ is used. This gives $(\ref{3.04})$.

 Computing the discrete inner product of $(\ref{2.02})$ with $\delta^{+}_{t}U^{n}$, then taking the real part, we obtain
\begin{equation}\label{3.06}
  -\frac{1}{2\tau}[(|U^{n+1}|^{2}_{H^{\alpha/2}_{h}}-\frac{\beta}{2}\|U^{n+1}\|^{4}_{l^{4}_{h}})-(|U^{n}|^{2}_{H^{\alpha/2}_{h}}-\frac{\beta}{2}\|U^{n}\|^{4}_{l^{4}_{h}})]
  =0,\quad t_{n}\in \Omega_{\tau},
\end{equation}
where Lemma $\ref{3.01}$ is used. This yields $(\ref{3.05})$.
\subsection{A priori estimate}
\begin{theorem}\label{3.12}
Then numerical solution of scheme $(\ref{2.02})$ is bounded in the following sense
\begin{equation}\label{3.07}
  \|U^{n}\|_{h}\leq C_{1},\quad \quad |U^{n}|_{H^{\alpha/2}_{h}}\leq C_{2},\quad\quad \|U^{n}\|_{l^{\infty}_{h}}\leq C_{3},\quad 0\leq n\leq N,
\end{equation}
where $C_{1},C_{2},C_{3}$ are some positive constants.
\end{theorem}
{\bf Proof.}\quad The proof is similar to that in~[\citealp[Theorem 3.2]{wang2016point}]. The mass conservation $(\ref{3.04})$ implies the first inequality in $(\ref{3.07})$ immediately if we choose $\|U^{0}\|_{h}\leq C_{1}.$\\
\indent Next, we prove the second inequality by the energy conservation $(\ref{3.05})$. If $\beta\leq 0$, according to the second term of $(\ref{3.02})$ and
energy conservation $(\ref{3.05})$, we can get the result straightforwardly. If $\beta>0$,
 In the view of Lemma $\ref{lemma3}$ with $\frac{1}{4}<\sigma_{0}<\frac{\alpha}{4}$, Young's inequality
and the first inequality in $(\ref{3.07})$, we obtain
\begin{equation}\label{3.08}
  \begin{aligned}
  \|U^{n}\|^{4}_{l^{4}_{h}}\leq C_{\sigma_{0}}\|U^{n}\|^{\frac{8\sigma_{0}}{\alpha}}_{H^{\alpha/2}_{h}}\|U^{n}\|^{4-\frac{8\sigma_{0}}{\alpha}}_{h}&\leq
  C_{\sigma_{0}}(\varepsilon|U^{n}|^{2}_{H^{\alpha/2}_{h}}+\varepsilon\|U^{n}\|^2_{h}+C(\varepsilon)),
  \end{aligned}
\end{equation}
where $\varepsilon$ is any arbitrary positive constant.   Combing the second term of $(\ref{3.02})$ with $(\ref{3.08})$, the energy conservation $(\ref{3.05})$ imply
 \begin{equation}\label{3.09}
  \begin{aligned}
  |U^{n}|^{2}_{H^{\alpha/2}_{h}}&= \frac{\beta}{2}\|U^{n}\|^{4}_{l^{4}_{h}}+E^{0}\\
  &\leq  \frac{\beta}{2}C_{\sigma_{0}}(\varepsilon|U^{n}|^{2}_{H^{\alpha/2}_{h}}+\varepsilon\|U^{n}\|^2_{h}+C(\varepsilon))
  +E^{0}.
  \end{aligned}
\end{equation}
Taking $\varepsilon=\frac{1}{ \beta C_{\sigma_{0}}}$, we have
\begin{equation}\label{3.10}
  |U^{n}|^{2}_{H^{\alpha/2}_{h}}\leq  \|U^{0}\|^{2}_{h}+\beta C_{\sigma_{0}}C(\varepsilon)+  2E^{0}:=C^{2}_{2},\quad 0\leq n\leq N_{t}.
\end{equation}
This implies the second inequality of $(\ref{3.02}).$\\
\indent Finally, combining the first two inequality in $(\ref{3.07})$ with Lemma $\ref{lemma1}$, we get the third inequality in $(\ref{3.07})$, that is,
\begin{equation}\label{3.11}
  \|U^{n}\|^{2}_{l^{\infty}_{h}}\leq C^{2}_{\sigma}(\|U^{n}\|^{2}_{h}+|U^{n}|^{2}_{H^{\alpha/2}_{h}})\leq C^{2}_{\sigma}(C^{2}_{1}+C^{2}_{2}):=C^{2}_{3},\quad 0\leq n\leq N_{t}.
\end{equation}
Thus the proof is completed.

\subsection{Existence}
\begin{theorem}
The nonlinear equation system in scheme $(\ref{2.02})$ is solvable.
\end{theorem}
{\bf Proof.}\quad The argument of the existence for the solution relies on the Browder fixed point theorem (see [\citealp[11]{akrivis1991fully}]). Here we omit the proof for brevity.
% In fact, since the discrete fractional Laplacian $(-\bigtriangleup)^{\frac{\alpha}{2}}_d$ is also self-adjoint~(see $(\ref{2.05})$), it is not difficult for us to apply the same method~(see Theorem 4.1 and Theorem 5.2 in  $\cite{wang2015energy}$) with lemma $\ref{2.04}$ and theorem $\ref{3.12}$ to obtain the proof of the uniquence. Here we omit the proof for brevity.

\section{Convergence of the scheme}
In this section, we will establish error estimate of $(\ref{2.02})$ in the discrete $L^{\infty}$ norm. For simplicity, we let $\Omega=[0,2\pi]$ and assume that $C^{\infty}_{p}(\Omega)$ is a set of infinitely differentiable functions with $2\pi$-period defined on $\Omega$. $H^{r}_{p}(\Omega)$ is the closure of $C^{\infty}_{p}(\Omega)$ in $H^{r}(\Omega)$. The semi-norm and the norm of $H^{r}_{p}(\Omega)$ are denoted by $|\cdot|$ and $\|\cdot\|_{r}$, respectively.  \\
\indent For the given even $N$, we introduce the projection space
$$S_{N}=\{u|u(x)=\sum_{|l|\leq N/2}\widehat{u}_{l}e^{\text{i}l x}\},$$
and the interpolation space
$$S^{''}_{N}=\{u|u(x)=\sum_{|l|\leq N/2}\frac{\widehat{u}_{l}}{c_{l}}\widehat{u}_{l}e^{\text{i}l x},\quad \widehat{u}_{-\frac{N}{2}}=\widehat{u}_{\frac{N}{2}},\}$$
where $c_{l}=1,|l|<\frac{N}{2},c_{-\frac{N}{2}}=c_{\frac{N}{2}}=2.$

 It is clear that $S^{''}_{N}\subseteq S_{N}.$ We denote by $P_{N}:L^{2}(\Omega)\rightarrow S_{N}$ as the orthogonal projection operator and recall the interpolation operator $I_{N}:L^{2}({\Omega})\rightarrow S^{''}_{N}.$  Further, $P_{N}$ and $I_{N}$ satisfy:
\begin{equation*}
  \begin{aligned}
 &1.P_{N}\partial_{x}u=\partial_{x}P_{N}u,~I_{N}\partial_{x}u\neq \partial_{x}I_{N}u.\\
 &2.P_{N}u=u,~\forall u=S_{N},I_{N}u=u,\forall u\in S^{''}_{N}.
  \end{aligned}
\end{equation*}
\begin{lemma}(\cite{gong2017conservative}) \label{lema}
 For $u\in S^{''}_{N},\|u\|\leq \|u\|_{h}\leq 2\|u\|.$
\end{lemma}
\begin{lemma} \label{lema1}
(\cite{canuto1982approximation}) If $0\leq l\leq r$ and $u\in H^{r}_{p}(\Omega),$ then
%\begin{equation}\label{lema1}
  \begin{gather}
    \|P_{N}u-u\|_{l}\leq CN^{l-r}|u|_{r},\\
    \|P_{N}u\|_{l}\leq C\|u\|_{l}.
  \end{gather}
 In addition, if $r>\frac{1}{2}$, then
 \begin{gather}
    \|I_{N}u-u\|_{l}\leq CN^{l-r}|u|_{r},\\
    \|I_{N}u\|_{l}\leq C\|u\|_{l}.
  \end{gather}

\end{lemma}
\begin{lemma} \label{lema2}
(\cite{gong2017conservative}) For $u\in H^{r}_{p}(\Omega),r>\frac{1}{2},$ let $u^{\ast}=P_{N-2}u,~$ then $\|u^{\ast}-u\|_{h}\leq CN^{-r}|u|_{r}.$
\end{lemma}

\begin{lemma} \label{lema3}
 For $u\in H^{r}_{p}(\Omega),r>\frac{1}{2},$ let $u^{\ast}=P_{N-2}u,~$ then
$|u^{\ast}-u|_{H^{\alpha/2}_{h}}\leq CN^{\alpha/2-r}|u|_{r}.$
\end{lemma}
{\bf Proof.}\quad According to Lemma $\ref{lemu}$, we have
\begin{equation}
  \begin{aligned}
|u^{\ast}-u|_{H^{\alpha/2}_{h}}&=((-\bigtriangleup)^{\frac{\alpha}{2}}_{d}(u^{\ast}-u),u^{\ast}-u)^{\frac{1}{2}}_{h}\\
&\leq \|(-\bigtriangleup)^{\frac{\alpha}{2}}_{d}(u^{\ast}-u)\|^{\frac{1}{2}}_{h}\|u^{\ast}-u\|^{\frac{1}{2}}_{h}\\
&=\|(-\bigtriangleup)^{\frac{\alpha}{2}}(I_{N}(u^{\ast}-u))\|^{\frac{1}{2}}_{h}\|u^{\ast}-u\|^{\frac{1}{2}}_{h}.
  \end{aligned}
\end{equation}
Together with Lemma $\ref{lema}$ and Lemma $\ref{lema1}$, we can deduce
\begin{equation}
  \begin{aligned}\label{a12}
\|(-\bigtriangleup)^{\frac{\alpha}{2}}(I_{N}(u^{\ast}-u))\|_{h}&=\|I_{N}[(-\bigtriangleup)^{\frac{\alpha}{2}}(I_{N}(u^{\ast}-u))]\|_{h}\\
&\leq \sqrt{2}\|I_{N}[(-\bigtriangleup)^{\frac{\alpha}{2}}(I_{N}(u^{\ast}-u))]\|\\
&\leq C\|(-\bigtriangleup)^{\frac{\alpha}{2}}(I_{N}(u^{\ast}-u))\|\\
&\leq C\|(I_{N}(u^{\ast}-u)\|_{\alpha}\\
&\leq C\|u^{\ast}-u\|_{\alpha}\leq CN^{\alpha-r}|u|_{r}.\\
  \end{aligned}
\end{equation}
Then, we can deduce from $(\ref{a12})$ and Lemma $\ref{lema2}$ that
$$|u^{\ast}-u|_{H^{\alpha/2}_{h}}\leq CN^{\alpha/2-r}|u|_{r}.$$
\begin{lemma} \label{lema4}
(\cite{sun2010convergence}) For time sequence $w= \{w^{0},w^{1},\cdots,w^{n}\},g=\{g^{\frac{1}{2}},g^{\frac{1}{2}},\cdots,g^{n-\frac{1}{2}}\}$ \\
\begin{equation}
  \begin{aligned}
|2\tau\sum_{l=1}^{n}g^{l-\frac{1}{2}}\cdot\delta^{+}_{t}w^{l}|&\leq (\tau\sum_{l=1}^{n-1}|\delta^{+}_{t}g^{l-\frac{1}{2}}|^{2}+\tau\sum_{l=1}^{n-1}|w^{l}|^{2})\\
&+|g^{n-\frac{1}{2}}|^{2}+|w^{n}|^{2}+|g^{\frac{1}{2}}|^{2}+|w^{0}|^{2}.
\end{aligned}
\end{equation}
\end{lemma}
\begin{lemma} \label{lema41}
(Gronwall Inequality(\cite{zhou1991applications})). Suppose that the discrete function $\{\omega^{n}|n=0,1,2\cdots N;$ $N\tau=T\}$ is nonnegative and satisfies the recurrence formula
\begin{equation}
  \begin{aligned}
\omega^{n}-\omega^{n-1}\leq A\omega^{n}\tau+B\omega^{n-1}\tau+C_{n}\tau,
\end{aligned}
\end{equation}
where $A,B$ and $C_{n}(n=1,2,\cdots)$ are nonnegative constants. Then
\begin{equation}
  \begin{aligned}
\max\limits_{0\leq n\leq N}|\omega^{n}|\leq (\omega^{0}+\sum_{k=1}^{N}C_{k}\tau)e^{2(A+B)T},
\end{aligned}
\end{equation}
where $\tau$ is satisfied $(A+B)\tau \leq \frac{N-1}{2N}(N>1)$.
\end{lemma}

\begin{lemma} \label{lema5}
(Gronwall Inequality\rm(\cite{zhou1991applications})). Suppose that the discrete function $\{\omega^{n}|n=0,1,2\cdots N;$ $N\tau=T\}$ is nonnegative and satisfies the inequality
\begin{equation}
  \begin{aligned}
\omega^{n}\leq A+\tau\sum_{l=1}^{n}B_{l}\omega^{l},
\end{aligned}
\end{equation}
where $A$ and $B_{l}(l=1,2,\cdots)$ are nonnegative constants. Then
\begin{equation}
  \begin{aligned}
\max\limits_{0\leq n\leq N}|\omega^{n}|\leq Ae^{2\sum_{l=1}^{N}B_{l}\tau},
\end{aligned}
\end{equation}
where $\tau$ is satisfied $\tau(\max\limits_{l=0,1,\cdots,N}B_{l})\leq \frac{1}{2}.$

\end{lemma}

\begin{lemma} \label{lema6}(\cite{sun2010convergence})
For any complex numbers $U,V,u,v,$ the following inequality holds
  \begin{align}
\|U|^{2}V-|u|^{2}v|\leq (\max\{|U|,|V|,|u|,|v|\})^{2}\cdot(2|U-u|+|V-v|).
\end{align}
\end{lemma}
\begin{theorem}\label{th1}
 We assume that the continuous solution $u$ of $(\ref{aa1})$ satisfies
\begin{equation}
  \begin{aligned}
u(x,t)\in C^{4}(0,t;H^{r}_{p}(\Omega)),\quad r>1,
\end{aligned}
\end{equation}
then the solution $U^{n}$ of $(\ref{2.02})$ is unconditionally convergent with order of $O(\tau^{2}+N^{\alpha/2-r})$ in the discrete $L^{\infty}$ norm.
\end{theorem}
{\bf Proof.}\quad We denote
\begin{equation}\label{4.01}
  \begin{aligned}
  u^{\ast}=P_{N-2}u,~~f=f(u)=\beta|u|^{2}u,\quad f^{\ast}=P_{N-2}f.
\end{aligned}
\end{equation}
The projection equation of $(\ref{aa1})$ is
\begin{equation}\label{4.02}
  \begin{aligned}
  \text{i}\partial_{t}u^{\ast}-(-\triangle)^{\frac{\alpha}{2}}u^{\ast}+f^{\ast}=0.
\end{aligned}
\end{equation}
We define
\begin{equation}\label{4.03}
  \begin{aligned}
  \xi^{n+\frac{1}{2}}_{j}=i\delta^{+}_{t}u^{\ast n}_{j}-(D_{\alpha}u^{\ast (n+1/2)})_{j}+f^{\ast}_{j}.
\end{aligned}
\end{equation}
Since $u^{\ast}\in S_{N},~~(-\Delta)^{\frac{\alpha}{2}}u^{\ast}(x_{j},t_{n})=(D_{\alpha}u^{\ast (n+1/2)})_{j}$, we obtain
\begin{equation}\label{4.04}
  \begin{aligned}
  \xi^{n+\frac{1}{2}}_{j}=\text{i}(\delta^{+}_{t}u^{\ast n}_{j}- \partial_{t}u^{\ast (n+1/2)}_{j}),
\end{aligned}
\end{equation}
and
\begin{equation}\label{4.40}
  \begin{aligned}
    \delta^{+}_{t}\xi^{n}_{j}&=\frac{\xi^{n+\frac{1}{2}}_{j}-\xi^{n-\frac{1}{2}}_{j}}{\tau}\\
    &=\frac{\text{i}}{\tau^{2}}(u^{\ast n+2}_{j}-2u^{\ast n+1}_{j}+u^{\ast n}_{j})-\frac{\text{i}}{2\tau}(\partial_{t}u^{\ast n+2}_{j}-\partial_{t}u^{\ast n}_{j}).
  \end{aligned}
\end{equation}
Using the Taylor expansion, we obtain
\begin{equation}\label{4.05}
  \begin{aligned}
  |\xi^{n+\frac{1}{2}}_{j}|\leq C\tau^{2},
\end{aligned}
\end{equation}
and
\begin{equation}\label{4.41}
  \begin{aligned}
  |\delta^{+}_{t}\xi^{n}_{j}|\leq C\tau^{2}.
\end{aligned}
\end{equation}
for some constant \emph C.

Denote $e^{n}_{j}=u^{\ast n}_{j}-U^{n}_{j}.$~Subtracting $(\ref{2.02})$ from $(\ref{4.02})$ yields the following error equation
  \begin{gather}\label{4.06}
  \xi^{n+\frac{1}{2}}=\text{i}\delta^{+}_{t}e^{n}-D_{\alpha}e^{n+1/2}+G^{n},\\
  e^{0}=u^{\ast 0}-u^{0},
\end{gather}
where
$$G^{n}_{j}=f^{\ast n+1/2}_{j}-F(U^{n}_{j},U^{n+1}_{j}).$$
Denoting
\begin{equation}\label{4.07}
  \begin{aligned}
 &(G_{1})^{n}_{j}=f^{\ast n+1/2}_{j}-f^{n+1/2}_{j},\quad (G_{2})^{n}_{j}=f^{n+1/2}_{j}-F(u^{n}_{j},u^{n+1}_{j}),\\
 &(G_{3})^{n}_{j}=F(u^{n}_{j},u^{n+1}_{j})-F(u^{\ast n}_{j},u^{\ast n+1}_{j}),\quad (G_{4})^{n}_{j}=F(u^{\ast n}_{j},u^{\ast n+1}_{j})-F(U^{n}_{j},U^{n+1}_{j}),
\end{aligned}
\end{equation}
we have
$G^{n}_{j}=(G_{1})^{n}_{j}+(G_{2})^{n}_{j}+(G_{3})^{n}_{j}+(G_{4})^{n}_{j}.$\\
According to Lemma $\ref{lema1}$, we have
\begin{equation}\label{4.08}
  \begin{aligned}
 \|G^{n}_{1}\|_{h}\leq CN^{-r}.
\end{aligned}
\end{equation}
By the Taylor expansion, we can see that
\begin{equation}\label{4.09}
  \begin{aligned}
 \|G^{n}_{2}\|_{h}\leq C\tau^{2}.
\end{aligned}
\end{equation}
From Lemma $\ref{lema6}$,~we deduced that
\begin{equation}\label{4.10}
  \begin{aligned}
 |(G_{3})^{n}_{j}|\leq C(|u^{n}_{j}-u^{\ast n}_{j}|+|u^{n+1}_{j}-u^{\ast n+1}_{j}|).
\end{aligned}
\end{equation}
This, together with Lemma $\ref{lema2}$ gives
\begin{equation}\label{4.11}
  \begin{aligned}
 \|G^{n}_{3}\|_{h}\leq CN^{-r}.
\end{aligned}
\end{equation}
From the definition of $G^{n}_{4}$, we have
\begin{equation}\label{4.12}
  \begin{aligned}
 (G_{4})^{n}_{j}&=\frac{\beta}{2}(|u^{\ast n}_{j}|^{2}+|u^{\ast n+1}_{j}|^{2})u^{\ast n+\frac{1}{2}}_{j}-\frac{\beta}{2}(|U^{n}_{j}|^{2}+|U^{ n+1}_{j}|^{2})U^{ n+\frac{1}{2}}_{j}\\
 &=\frac{\beta}{2}(|u^{\ast n}_{j}|^{2}+|u^{\ast n+1}_{j}|^{2}-|U^{n}_{j}|^{2}-|U^{n+1}_{j}|^{2})u^{\ast n+\frac{1}{2}}_{j}+\frac{\beta}{2}(|U^{n}_{j}|^{2}+|U^{n+1}_{j}|^{2})e^{n+\frac{1}{2}}_{j}\\
 &:=(G_{41})^{n}_{j}+(G_{42})^{n}_{j},
\end{aligned}
\end{equation}
where
\begin{equation}\label{4.13}
  \begin{aligned}
&(G_{41})^{n}_{j}=\frac{\beta}{2}(|u^{\ast n}_{j}|^{2}+|u^{\ast n+1}_{j}|^{2}-|u^{\ast n}_{j}-e^{n}_{j}|^{2}-|u^{\ast n+1}_{j}-e^{n+1}_{j}|^{2})u^{\ast n+\frac{1}{2}}_{j}\\
&\quad ~~~\quad=\frac{\beta}{2}(u^{\ast n}_{j}\overline{e}^{n}_{j}+\overline{u^{\ast}}^{n}_{j}e^{n}_{j}+u^{\ast n+1}_{j}\overline{e}^{n+1}_{j}+\overline{u^{\ast}}^{n+1}_{j}e^{n+1}_{j}-|e^{n}_{j}|^{2}-|e^{n+1}_{j}|^{2})u^{\ast n+1/2}_{j}.\\
&(G_{42})^{n}_{j}=\frac{\beta}{2}(|U^{n}_{j}|^{2}+|U^{n+1}_{j}|^{2})e^{n+\frac{1}{2}}_{j}.
\end{aligned}
\end{equation}
Computing the discrete inner product of $(\ref{4.06})$ with $e^{n+1/2}$, then taking the imaginary part, we obtain
\begin{equation}\label{4.14}
  \begin{aligned}
\frac{1}{2\tau}(\|e^{n+1}\|^{2}_{h}-\|e^{n}\|^{2}_{h})+\textrm{Im}(G^{n}_{1}+G^{n}_{2}+G^{n}_{3}+G^{n}_{4},e^{n+1/2})_{h}=\textrm{Im}(\xi^{n+\frac{1}{2}},e^{n+1/2})_{h},
\end{aligned}
\end{equation}
Using Cauchy-Schwartz inequality, we obtain
\begin{gather}\label{4.17}
|(G^{n}_{s},e^{n+1/2})_{h}|\leq \frac{1}{2}\|G^{n}_{s}\|^{2}_{h}+\frac{1}{4}(\|e^{n}\|^{2}_{h}+\|e^{n+1}\|^{2}_{h}),~s=1,2,3,\\
|\textrm{Im}(G^{n}_{4},e^{n+1/2})_{h}|
=|\textrm{Im}(G^{n}_{41},e^{n+1/2})_{h}|\leq C(\|e^{n}\|^{2}_{h}+\|e^{n+1}\|^{2}_{h}+\|e^{n}\|^{4}_{l^{4}_{h}}+\|e^{n+1}\|^{4}_{l^{4}_{h}}),\label{4.18}\\
|(\xi^{n+\frac{1}{2}},e^{n+1/2})_{h}|\leq \frac{1}{2}|\xi^{n+\frac{1}{2}}|^{2}_{h}+\frac{1}{4}(\|e^{n}\|^{2}_{h}+\|e^{n+1}\|^{2}_{h}).\label{4.21}
\end{gather}
According to theorem $\ref{3.12}$, Lemma~$\ref{lema}$ and Lemma~$\ref{lema1}$, we can get
\begin{gather}\label{4.15}
\|e^{n}\|_{h}=\|u^{\ast n}-U^{n}\|_{h}\leq \|u^{\ast}\|_{h}+\|U^{n}\|_{h}\leq \sqrt{2}\|u^{\ast n}\|+\|U^{n}\|_{h}\leq C.\\
|e^{n}|_{H^{\alpha/2}_{h}}=|u^{\ast n}-U^{n}|_{H^{\alpha/2}_{h}}\leq |u^{\ast n}|_{H^{\alpha/2}_{h}}+|U^{n}|_{H^{\alpha/2}_{h}}\leq C \label{4.19}.
\end{gather}
From $(\ref{4.15})\verb|-|(\ref{4.19})$ and Lemma $\ref{lemma3}$, we have
\begin{equation}\label{4.16}
  \begin{aligned}
\|e^{n}\|^{4}_{l^{4}_{h}}\leq C\|e^{n}\|^{2}_{h}.
\end{aligned}
\end{equation}
This, together with $(\ref{4.18})$, gives
\begin{align}\label{4.20}
 & |\textrm{Im}(G^{n}_{4},e^{n+1/2})_{h}|\leq C(\|e^{n}\|^{2}_{h}+\|e^{n+1}\|_{h}),\\
 & \|G^{n}_{4}\|^{2}_{h}\leq C(\|e^{n}\|^{2}_{h}+\|e^{n+1}\|_{h}).
\end{align}
Substituting $(\ref{4.17}),(\ref{4.21}),(\ref{4.20})$ into $(\ref{4.14})$ yields
\begin{equation}\label{4.22}
  \frac{1}{2\tau}(\|e^{n+1}\|^{2}_{h}-\|e^{n}\|^{2}_{h})\leq C(\|e^{n}\|^{2}_{h}+\|e^{n+1}\|^{2}_{h})+\frac{1}{2}(\|G^{n}_{1}\|^{2}_{h}+\|G^{n}_{2}\|^{2}_{h}+\|G^{n}_{3}\|^{2}_{h}+\|\xi^{n+\frac{1}{2}}\|^{2}_{h}).
\end{equation}
This together with Lemma $\ref{lema41}$ and $(\ref{4.08}),(\ref{4.09}),(\ref{4.11})$, gives that, for a sufficiently small $\tau$,
\begin{equation}\label{4.23}
  \|e^{n}\|^{2}_{h}\leq (\|e^{0}\|^{2}_{h}+CT(N^{-2r}+\tau^{4}))e^{4CT}.
\end{equation}
This together with
\begin{equation}\label{4.24}
  \|e^{0}\|_{h}=\|u^{\ast 0}-u^{0}\|_{h}\leq CN^{-r},
\end{equation}
gives
\begin{equation}\label{4.25}
  \|e^{n}\|_{h}\leq C(N^{-r}+\tau^{2}).
\end{equation}
With Lemma $\ref{lema2}$ and $(\ref{4.25})$, we can get
\begin{equation}\label{4.26}
  \|u^{n}-U^{n}\|_{h}\leq \|u^{n}-u^{\ast n}\|_{h}+\|u^{\ast n}-U^{n}\|_{h}\leq C(N^{-r}+\tau^{2}).
\end{equation}
Computing the discrete inner product of $(\ref{4.06})$ with $\delta^{+}_{t}e^{n}$, and taking the real part, we obtain
\begin{equation}\label{4.27}
  \begin{aligned}
\frac{1}{\tau}(|e^{n+1}|^{2}_{H^{\alpha/2}_{h}}-|e^{n}|^{2}_{H^{\alpha/2}_{h}})=-2\textrm{Re}(G^{n},\delta^{+}_{t}e^{n})_{h}+2\textrm{Re}(\xi^{n+\frac{1}{2}},\delta^{+}_{t}e^{n})_{h},
\end{aligned}
\end{equation}
where $(\ref{a8})$ is used. From $(\ref{4.06})$, we have
\begin{equation}\label{4.28}
  \delta^{+}_{t}e^{n}=-\text{i}D_{\alpha}e^{n+1/2}+\text{i}G^{n}_i-\text{i}\xi^{n+\frac{1}{2}}.
\end{equation}
Substituting $(\ref{4.28})$ into the first term on the right side of $(\ref{4.27})$, we have
\begin{equation}\label{4.29}
  \begin{aligned}
  2\textrm{Re}(G^{n},\delta^{+}_{t}e^{n})&=2\textrm{Re}(G^{n},-\text{i}D_{\alpha}e^{n+1/2}+\text{i}G^{n}-\text{i}\xi^{n+\frac{1}{2}})_{h}\\
  &=2\textrm{Im}(G^{n},D_{\alpha}e^{n+1/2})_{h}+2\textrm{Im}(G^{n},\xi^{n+\frac{1}{2}})_{h}.
\end{aligned}
\end{equation}
By virtue of $(\ref{2.07})$, we have
\begin{equation}\label{4.30}
  \begin{aligned}
  2\textrm{Im}(G^{n},D_{\alpha}e^{n+1/2})_{h}&\leq 2|G^{n}|_{H^{\alpha/2}_{h}}|e^{n+1/2}|_{H^{\alpha/2}_{h}}\\
  &\leq |G^{n}|^{2}_{H^{\alpha/2}_{h}}+|e^{n+1/2}|^{2}_{H^{\alpha/2}_{h}},
\end{aligned}
\end{equation}
and
\begin{equation}\label{4.31}
  2\textrm{Im}(G^{n},\xi^{n+\frac{1}{2}})_{h}\leq 2\|G^{n}\|_{h}\|\xi^{n+\frac{1}{2}}\|_{h}\leq \|G^{n}\|^{2}_{h}+\|\xi^{n+\frac{1}{2}}\|^{2}_{h}.
\end{equation}
The estimate of $|G^{n}|_{H^{\alpha/2}_{h}}$ is established as follows. According to Lemma $\ref{lema3}$, we have
\begin{equation}\label{4.32}
  |G^{n}_{1}|_{H^{\alpha/2}_{h}}\leq CN^{\alpha/2-r}.
\end{equation}
By Taylor formula, we can see that
\begin{equation}\label{4.33}
  |G^{n}_{2}|_{H^{\alpha/2}_{h}}\leq C\tau^{2}.
\end{equation}
With the similar argument for $G^{n}_{3}$ and $G^{n}_{4}$ as above, we can deduce that
\begin{equation}\label{4.34}
  |G^{n}_{3}|_{H^{\alpha/2}_{h}}\leq CN^{\alpha/2-r},
\end{equation}
and
\begin{equation}\label{4.35}
  \begin{aligned}
 |G^{n}_{4}|^{2}_{H^{\alpha/2}_{h}}\leq C(\|e^{n}\|^{2}_{h}+\|e^{n+1}\|^{2}_{h}+|e^{n}|^{2}_{H^{\alpha/2}_{h}}+|e^{n+1}|^{2}_{H^{\alpha/2}_{h}}),
\end{aligned}
\end{equation}
where Lemma $\ref{lema3}$ is used. Therefore, it follows from $(\ref{4.32})$ and $(\ref{4.35})$ that
\begin{equation}\label{4.36}
\begin{aligned}
  |G^{n}|^{2}_{H^{\alpha/2}_{h}}&\leq C(\|e^{n}\|^{2}_{h}+\|e^{n+1}\|^{2}_{h}+|e^{n}|^{2}_{H^{\alpha/2}_{h}}+|e^{n+1}|^{2}_{H^{\alpha/2}_{h}})\\
  &+C(N^{\alpha-2r}+\tau^{4}).
  \end{aligned}
\end{equation}
Hence we get
\begin{equation}\label{4.37}
\begin{aligned}
  2\textrm{Re}(G^{n},\delta^{+}_{t}e^{n})&\leq C(\|e^{n}\|^{2}_{h}+\|e^{n+1}\|^{2}_{h}+|e^{n}|^{2}_{H^{\alpha/2}_{h}}+|e^{n+1}|^{2}_{H^{\alpha/2}_{h}})\\
  &+C(N^{\alpha-2r}+\tau^{4}).
  \end{aligned}
\end{equation}
Substituting $(\ref{4.37})$  into $(\ref{4.27})$, we have
\begin{equation}\label{4.38}
\begin{aligned}
|e^{n+1}|^{2}_{H^{\alpha/2}_{h}}-|e^{n}|^{2}_{H^{\alpha/2}_{h}}&\leq C\tau(\|e^{n}\|^{2}_{h}+\|e^{n+1}\|^{2}_{h}+|e^{n}|^{2}_{H^{\alpha/2}_{h}}+|e^{n+1}|^{2}_{H^{\alpha/2}_{h}})\\
&+2\tau \textrm{Re}(\xi^{n+\frac{1}{2}},\delta^{+}_{t}e^{n})_{h}+C\tau(N^{\alpha-2r}+\tau^{4}).
  \end{aligned}
\end{equation}
Summing up  the superscript $n$ from $0$ to $M$ and then replacing $M$ by $n$, we have
\begin{equation}\label{4.39}
\begin{aligned}
|e^{n+1}|^{2}_{H^{\alpha/2}_{h}}-|e^{0}|^{2}_{H^{\alpha/2}_{h}}&\leq \tau\sum_{l=0}^{n}C(\|e^{l}\|^{2}_{h}+\|e^{l+1}\|^{2}_{h}+|e^{l}|^{2}_{H^{\alpha/2}_{h}}+|e^{l+1}|^{2}_{H^{\alpha/2}_{h}})\\
&+2\tau \sum_{l=0}^{n}\textrm{Re}(\xi^{l+\frac{1}{2}},\delta^{+}_{t}e^{l})_{h}+CT(N^{\alpha-2r}+\tau^{4}).
  \end{aligned}
\end{equation}
With the Lemma $\ref{lema4},(\ref{4.05})$ and $(\ref{4.41})$, we have
\begin{equation}\label{4.42}
\begin{aligned}
|2\tau \sum_{l=0}^{n}\textrm{Re}(\xi^{l+\frac{1}{2}},\delta^{+}_{t}e^{l})_{h}|&=|\textrm{Re}(2\tau\sum_{l=0}^{n}h\sum_{j=0}^{N-1}\xi^{l+\frac{1}{2}}_{j}(\delta^{+}_{t}e^{l}_{j}))|\\
&\leq (h\sum_{j=0}^{N-1}2\tau\sum_{l=1}^{n}\xi^{l+\frac{1}{2}}(\delta^{+}_{t}e^{l}_{j}))\\
&\leq (\tau\sum_{l=1}^{n-1}\|\delta^{+}_{t}\xi^{l-\frac{1}{2}}\|_{h}+\tau\sum_{l=0}^{n}\|e^{l}\|^{2}_{h})\\
&+\|\xi^{n+\frac{1}{2}}\|^{2}_{h}+\|e^{n+1}\|^{2}_{h}+\|\xi^{\frac{1}{2}}\|^{2}_{h}+\|e^{0}\|^{2}_{h}\\
&\leq C(N^{-r}+\tau^{4}).
  \end{aligned}
\end{equation}
By virtue of $(\ref{4.26})$, we deduce from $(\ref{4.39})$ and $(\ref{4.42})$
\begin{equation}\label{4.43}
\begin{aligned}
|e^{n+1}|^{2}_{H^{\alpha/2}_{h}}&\leq |e^{0}|^{2}_{H^{\alpha/2}_{h}}+\tau\sum_{l=0}^{n}C(|e^{l}|^{2}_{H^{\alpha/2}_{h}}+|e^{l+1}|^{2}_{H^{\alpha/2}_{h}})+C(N^{\alpha-2r}+\tau^{4})\\
&\leq \tau\sum_{l=1}^{n+1}C|e^{l}|^{2}_{H^{\alpha/2}_{h}}+C(N^{\alpha-2r}+\tau^{4}),
  \end{aligned}
\end{equation}
where $|e^{0}|_{H^{\alpha/2}_{h}}\leq CN^{\alpha-2r}$ is used.\\
\indent Applying Lemma \ref{lema5} to $(\ref{4.43})$, we get
\begin{equation}\label{4.44}
  |e^{n}|^{2}_{H^{\alpha/2}_{h}}\leq (N^{\alpha-2r}+\tau^{4})e^{2\sum_{k=1}^{N}C\tau},
\end{equation}
where $\tau$ is sufficiently small, such that $C\tau\leq \frac{1}{2}$. With Lemma $\ref{lema3}$ and $(\ref{4.44})$, we can prove
\begin{equation}\label{4.45}
  |u^{n}-U^{n}|^{2}_{H^{\alpha/2}_{h}}\leq |u-u^{\ast}|_{H^{\alpha/2}_{h}}+|u^{\ast}-U|_{H^{\alpha/2}_{h}}\leq C(N^{\alpha-2r}+\tau^{4}).
\end{equation}
Finally, thanks to Lemma \ref{lemma1}, we obtain
\begin{equation}\label{4.46}
  \|u^{n}-U^{n}\|_{l^{\infty}_{h}}\leq C(N^{\alpha/2-r}+\tau^{2}).
\end{equation}
This completes the proof.
\section{Numerical examples}
In this section, we test the numerical accuracy and discrete conservation laws of the fully discrete pseudo-spectral scheme $(\ref{2.02})$. Similar to that in ~\cite{gong2017conservative}, we can use the fixed point iteration method and the fast Fourier transform~(FFT) to solve the nonlinear system defined in scheme  $(\ref{2.02})$.
% for each time step $n$, the fixed point iteration method is applied to solve the nonlinear system defined in scheme  $(\ref{2.02})$, that is 
%\begin{align}
%\text{i}\frac{U^{n+1}_{s+1}-U^{n}}{\tau}-D_{\alpha}\frac{U^{n+1}_{s+1}+U^{n}}{2}+F(U^{n},U^{n+1}_{s})=0,
%\end{align}
%where the initial iteration $U^{n+1}_{0}=U^{n}.$ Moreover, the iterative solution satisfies
%\begin{align}
%	\|U^{n+1}_{s+1}-U^{n+1}_{s}\|_{h,\infty}\leq 10^{-14}.
%\end{align}
%Since $D_{\alpha}U^{n+1}_{s+1}=\mathcal{F}^{-1}_{d}\bm{\wedge_{\alpha}}\mathcal{F}_{d}U^{n+1}_{s+1}$,
% we can use the fast Fourier transform~(FFT) to solve the linear equations efficiently for a fixed iteration step s. 
For the convergence rate, we use the formula
\begin{align}
 Order=\frac{\text{ln}(error_{1}/error_{2})}{\text{ln}(\tau_{1}/\tau_{2})},
\end{align}
where $\tau_{l},$ ${ error}_{l},(l=1,2)$ are step sizes and corresponding errors, respectively. The relative errors of energy and mass are defined as
\begin{align}
RH^{n}=|(H^{n}-H^{0})/(H^{0})|,~~RM^{n}=|(M^{n}-M^{0})/(H^{0})|.
\end{align}
%where $H^{n}$ and $M^{n}$ denote the energy and mass at $t=n\tau$, respectively.

\subsection{Example 1} Consider the FNLS equation~$(\ref{aa1})$ with plane wave solution \cite{duo2016mass}
\begin{align}
 u(x,t)=A{\rm exp}(\text{i}(\lambda x-\omega t)),~~~~\text{with} ~~\omega=|\lambda|^{\alpha}-\beta|A|^{2},
\end{align}
where the problem is solved on domain $[-\pi,\pi]$ with $\lambda=4, A=1,\beta=-2.$ 

Firstly, we test the accuracy and efficiency of numerical scheme $(\ref{2.02})$ for different power $\alpha$.   \tablename{~\ref{tab1}} indicates that the method is of second-order in time. \tablename{~\ref{tab2}} shows that the  spatial error is very small and almost negligible, and the error is dominated by the time discretization error. It confirms that, for sufficiently smooth problem, the Fourier pseudo-spectral method is of arbitrary order of accuracy. These numerical results confirm the accuracy of the numerical scheme $(\ref{2.02})$ in Theorem~\ref{th1}.

\begin{figure}[htbp]
	\centering
	\subfigure[Relative errors of energy]{
		\begin{minipage}[t]{0.5\linewidth}
			\centering
			\includegraphics[scale=0.5]{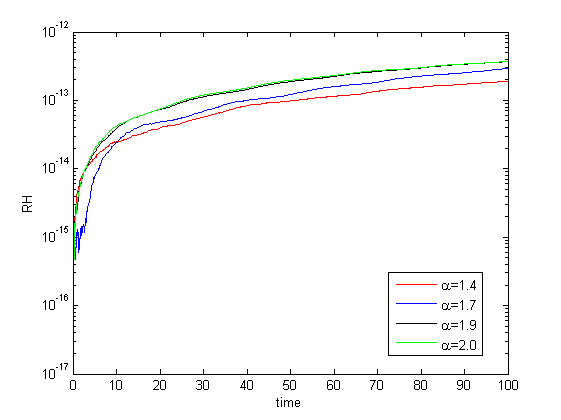}
			%\caption{fig1}
		\end{minipage}%
	}%
	\subfigure[Relative errors of mass]{
		\begin{minipage}[t]{0.5\linewidth}
			\centering
			\includegraphics[scale=0.5]{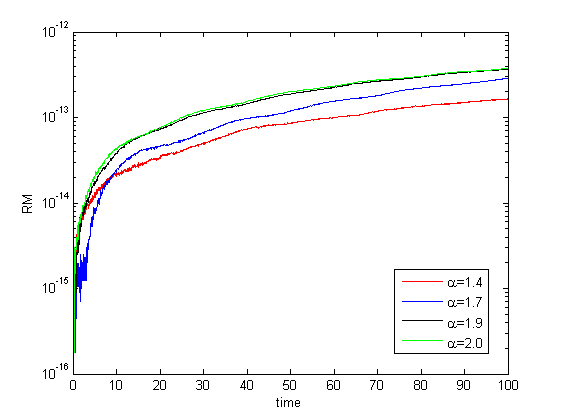}
			%\caption{fig2}
		\end{minipage}%
	}%
\caption{Relative errors of energy and mass with $N=32$, $\tau$=0.05 for different $\alpha$.}\label{fig1}%图的名字可以label
\end{figure}

Secondly, we testify the discrete conservation laws. Considering the iteration error and the slow-varying process of the solution, it is easy to select large-step $\tau=0.05$ in relative errors test of energy and mass. In Figure{~\ref{fig1}}, we depict the relative errors of energy and mass in a longer time interval. It is observed that the scheme $(\ref{2.02})$ conserves both energy and mass very well.
\begin{table}[htbp]
	\centering
	\caption{Convergence test in time for different $\alpha$ with $N=256$ at $T=1$.}\label{tab1}
	\setlength{\tabcolsep}{6mm}
	\begin{tabular}{ccccccc}
		\toprule[1pt]
		\rowcolor[gray]{0.9}  $\alpha$& $\tau$ & $L^{\infty}$&Rate&$L^{2}$& Rate \\
		1.4&0.05&   0.1529&-&  0.3831&-	\\
		&0.025& 0.0382&2.0017  &  0.0957&  2.0017\\
		&0.0125&0.0095&2.0095& 0.0238&2.0095\\ 
		&0.00625& 0.0024&  2.0069& 0.0059& 2.0069\\
		1.7&0.05 &0.4062&-& 1.0183&-\\
		&0.025&0.1041& 1.9644& 0.2609& 1.9644\\
		& 0.0125&0.0260&2.0013&0.0652&2.0013\\
		&0.00625 & 0.0065&  2.0049&  0.0162& 2.0049\\
		1.9&0.05& 0.7866&-& 1.9718&-\\
		&0.025&0.2104& 1.9027&0.5273&1.9027\\
		&0.0125&0.0530&1.9900&0.1328&1.9900\\
		&0.00625& 0.0132& 2.0024& 0.0331& 2.0024\\
		2.0&0.05&1.0794&-&2.7056&-\\
		&0.025&0.3011&1.8420&0.7547&1.8420\\
		&0.0125&0.0763&1.9805&0.1912&1.9806\\
		&0.00625&0.0191&2.0004&0.0478&2.0004\\
		\hline
	\end{tabular}
\end{table}
\begin{table}[htbp]
	\centering
	\caption{Convergence test in space for different $\alpha$ with $\tau$=1.0e-6 at $T=1$.}\label{tab2}
	\setlength{\tabcolsep}{12mm}
	\begin{tabular}{ccccccc}
		\toprule[1pt]
		\rowcolor[gray]{0.9}  $\alpha$& N & $L^{\infty}$&$L^{2}$  \\
		1.4& 16& 1.5612e-10&3.5615e-10\\
		& 32&3.0115e-10& 7.3557e-10\\
		&64&1.0579e-10&2.0756e-10\\
		&128& 1.0232e-10& 2.0650e-10 \\
		%&256&2.1464e-10&5.1767e-10\\
		1.7&16&1.2158e-10&2.9206e-10\\
		&32& 9.9284e-11&2.4356e-10\\
		&64& 2.0598e-10&4.4501e-10\\
		&128&1.7674e-10&4.0082e-10\\
		%&256&7.2828e-11& 1.7630e-10\\
		1.9&16&3.1580e-10& 7.8269e-10\\
		&32& 1.9745e-10&4.8347e-10\\
		&64&4.0320e-10& 9.4242e-10\\
		&128& 3.7175e-10&8.9973e-10\\
		%&256& 2.5302e-10&6.3172e-10\\
		2.0&16& 6.0690e-10&1.5007e-09\\
		&32& 1.9745e-10&4.8347e-10\\
		&64& 4.0320e-10& 9.4242e-10\\
		&128&3.7175e-10&8.9973e-10\\
		%&256& 2.5302e-10& 6.3172e-10\\
		\hline
	\end{tabular}
\end{table}
\begin{figure}[htbp]
	\centering
	\subfigure[Relative errors of energy]{
		\begin{minipage}[t]{0.5\linewidth}
			\centering
			\includegraphics[scale=0.5]{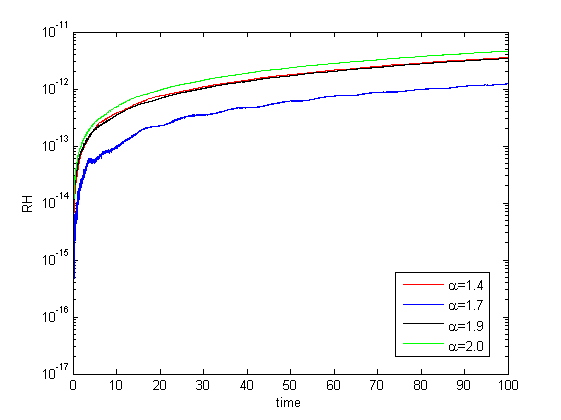}
			%\caption{fig1}
		\end{minipage}%
	}%
	\subfigure[Relative errors of mass]{
		\begin{minipage}[t]{0.5\linewidth}
			\centering
			\includegraphics[scale=0.5]{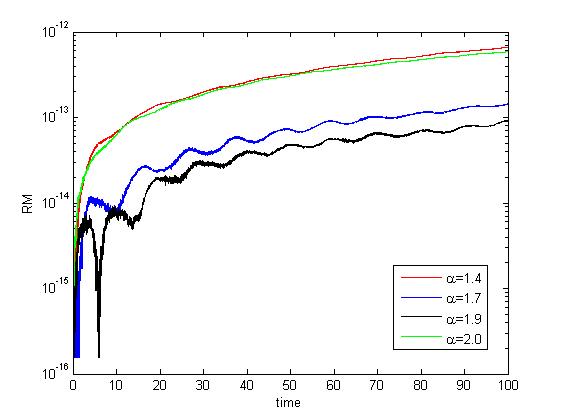}
			%\caption{fig2}
		\end{minipage}%
	}%
	\caption{Relative errors of energy and mass with $N=320$, $\tau$=0.01 for different $\alpha$.}\label{fig2}%图的名字可以label
\end{figure}

\subsection{Example 2} 
In this example, we study the dynamics of solitions in the one-dimensional FNLS equation
\begin{align}
	 \text{i}u_{t}-(-\Delta)^{\frac{\alpha}{2}} u+|u|^{2}u=0,
\end{align}
with initial condition
\begin{align}
u(x,0)={\rm sech}(\sqrt{2}x/2){\rm exp}(\text{i}x/2),~~~~x\in  [-20,20].
\end{align}
For $\alpha=2$, the exact solution is the solition on the whole real axis. In our simulations, we choose the spatial grid size $N=320$.

\begin{figure}[htbp]
	\centering
	\subfigure[$\alpha$=2, $\tau$=0.01]{
		\begin{minipage}[t]{0.5\linewidth}
			\centering
			\includegraphics[scale=0.5]{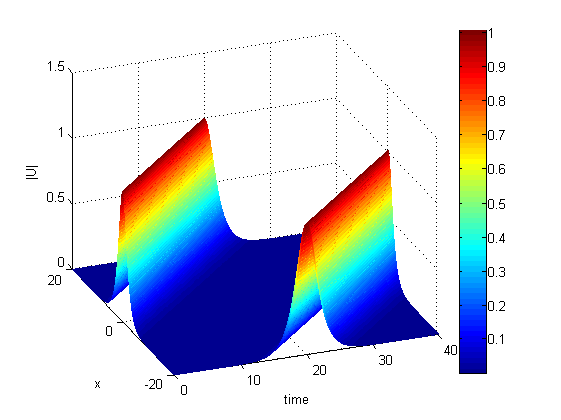}
			%\caption{fig1}
		\end{minipage}%
	}%
	\subfigure[$\alpha$=1.95, $\tau$=0.01]{
		\begin{minipage}[t]{0.5\linewidth}
			\centering
			\includegraphics[scale=0.5]{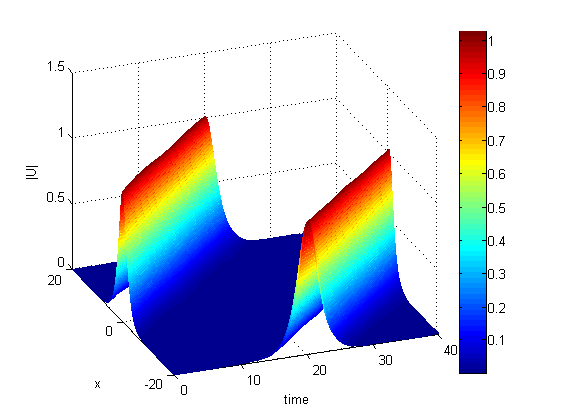}
			%\caption{fig2}
		\end{minipage}%
	}\\%
	\subfigure[$\alpha$=1.9, $\tau$=0.01]{
	\begin{minipage}[t]{0.5\linewidth}
		\centering
		\includegraphics[scale=0.5]{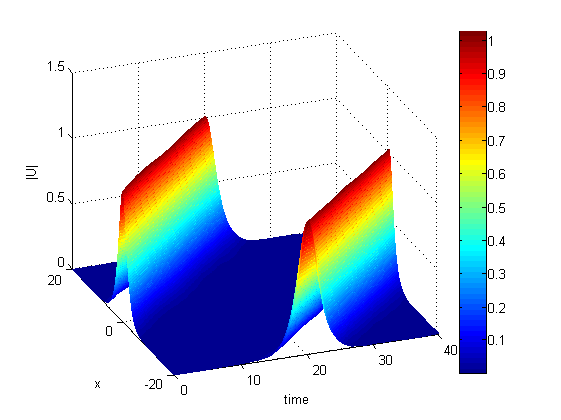}
		%\caption{fig2}
	\end{minipage}%
}%
	\subfigure[$\alpha$=1.7, $\tau$=0.01]{
	\begin{minipage}[t]{0.5\linewidth}
		\centering
		\includegraphics[scale=0.5]{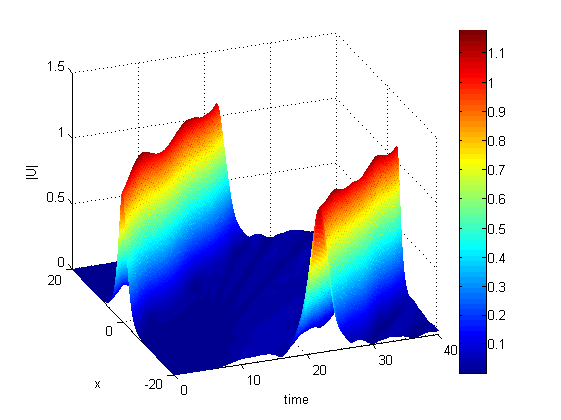}
		%\caption{fig2}
	\end{minipage}%
}\\%
	\subfigure[$\alpha$=1.5, $\tau$=0.01]{
	\begin{minipage}[t]{0.5\linewidth}
		\centering
		\includegraphics[scale=0.5]{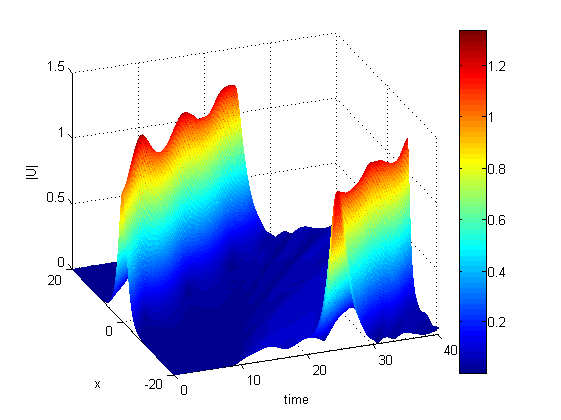}
		%\caption{fig2}
	\end{minipage}%
}%
	\subfigure[$\alpha$=1.3, $\tau$=0.01]{
	\begin{minipage}[t]{0.5\linewidth}
		\centering
		\includegraphics[scale=0.5]{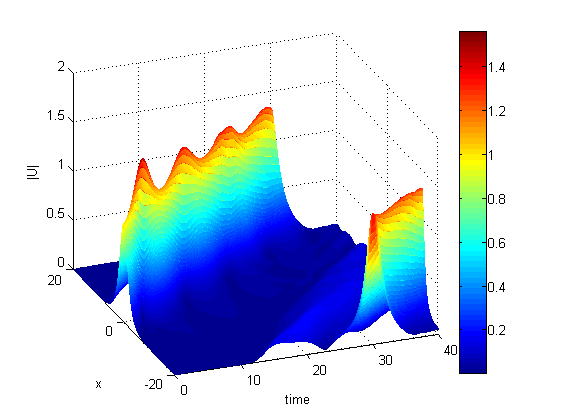}
		%\caption{fig2}
	\end{minipage}%
}%

	\caption{Evolution of the solitons with $N=320$, $\tau$=0.01 for different $\alpha$.}\label{fig3}%图的名字可以label
\end{figure}

In  Figure{~\ref{fig2}}, the relative errors of energy and mass are plotted, which show that energy and mass are conserved very well. Figure{~\ref{fig3}} presents the time evolution of the density $|u(x,t)|$ of the FNLS with different power $\alpha$. In the classical NLS with $\alpha=2$, the shape and velocity of soliton solutions are unchanged. In the FNLS, the shape of soliton solutions presents a slow changing process.
 When $\alpha$ decreases, the propagation of waves in the time-axis direction slows down with the elapse of time.
  We can also see that the solution of the FNLS behaves more like a wave with effects that might be described as
  \textquotedblleft
   interference\textquotedblright~arising from the long-range interactions of the fractional Laplacian. All these properties of FNLS may better simulate the shape of waves in physics.

\section{Conclusions}

In this paper,  unconditional convergence analysis of a conservative Fourier pseudo-spectral method for solving the FNLS equation is established. We introduce the discrete fractional Sobolev space $H^{\alpha/2}_{h}$ with a new discrete fractional Sobolev norm for the first time, and we also prove some lemmas for the new Sobolev norm. Based on these lemmas and energy method, a priori error estimate for the method can be estimated.
Then, we can prove that the conservative Fourier pseudo-spectral method is unconditionally convergent with order of $O(\tau^{2}+N^{\alpha/2-r})$ in the discrete $L^{\infty}$ norm.
  In fact, here we adopt a more direct analysis method so that unconditionally convergent results can be obtained without establishing semi-norm equivalence in error analysis.
  Furthermore, The method of error analysis in the discrete $L^{\infty}$ for the presented conservative Fourier pseudo-spectral
  scheme  can be extended to other PDEs involving fractional Laplacian, for example, the space fractional Allen-Cahn equation~\cite{song2016fractional} and the Klein-Gordon-Schr\"{o}dinger equation~\cite{huang2016global}.
  
\section*{Acknowledgments}

This work is supported by the National Key Research and Development Project of China (Grant No. 2017YFC0601505, 2018YFC1504205),
the National Natural Science Foundation of China (Grant No. 11771213, 61872422, 11971242, 11901513),
the Major Projects of Natural Sciences of University in Jiangsu Province of China (Grant No.18KJA110003),
and the Priority Academic Program Development of Jiangsu Higher Education Institutions.

\bibliography{main1}
\end{document}